\documentclass[12pt]{article}
\RequirePackage{etex} 
\usepackage[english]{babel}
\usepackage{graphicx}
\usepackage{comment}
\usepackage{amsmath,amsthm,amssymb}
\usepackage{times}
\usepackage{enumerate}
\usepackage{amsmath}
\usepackage{amsfonts}
\usepackage{pdfsync}
\usepackage{cleveref}
\usepackage{autonum}
\usepackage{enumitem}
\usepackage[usenames, dvipsnames]{color}

\newcommand{\R}{{\mathbb R}}
\newcommand{\N}{{\mathbb N}}

\DeclareMathOperator{\sgn}{sgn}

\DeclareMathOperator{\Lip}{Lip}
\def\titlerunning#1{\gdef\titrun{#1}}
\makeatletter
\def\author#1{\gdef\autrun{\def\and{\unskip, }#1}\gdef\@author{#1}}
\def\address#1{{\def\and{\\}\renewcommand{\thefootnote}{}%
\footnote {#1}}%
\markboth{\autrun}{\titrun}}
\makeatother
\def\email#1{e-mail: #1}
\def\subjclass#1{{\renewcommand{\thefootnote}{}%
\footnote{\emph{Mathematics Subject Classification (2010):} #1}}}

\theoremstyle{definition}
\newtheorem{theorem}{Theorem}[section]
\newtheorem{corollary}[theorem]{Corollary}
\newtheorem{lemma}[theorem]{Lemma}

\newtheorem{proposition}[theorem]{Proposition}
\newtheorem*{theorem*}{Theorem}
\usepackage{amsthm}

\theoremstyle{definition}
\newtheorem{definition}[theorem]{Definition}
\newtheorem{remark}[theorem]{Remark}

\newtheorem*{hypothesisH1'}{Hypothesis $\mathbf{(H_2')}$}
\newtheorem*{hypothesisH1''}{Hypothesis ($\mathcal S^{\infty}_{x_*}$)}
\newtheorem*{hypothesisS'}{Hypothesis ($\mathcal S^{\infty}_{x_*}$)}
\newtheorem*{hypothesisH'}{Hypothesis ($\mathcal A^{\infty}_{x_*}$)}

\newtheorem*{CondRy}{Condition (R$_y$)}
\newtheorem*{CondR}{Condition (R)}

\newtheorem*{problemP'}{Problem \textbf{(P$'$)}}

\numberwithin{equation}{section}

\begin{document}
\titlerunning{Occurrence of gap}
\title{Occurrence of gap for one-dimensional scalar autonomous functionals with one end point condition}
\author{
Rapha\"el  Cerf
\and
 Carlo Mariconda
\thanks{
This research was partially supported by a grant from  Padua University, grant SID 2018 ``Controllability, stabilizability and infimum gaps for control systems,'' prot. BIRD 187147, and was accomplished within the UMI Group TAA ``Approximation Theory and Applications.''  {C.M. wishes to thank Rapha\"el Cerf for the hospitality and invitation at the \'Ecole Normale Supérieure of Paris in May 2022, during the preparation of the paper and for the collaboration, more than 30 years after our last joint paper.}}}
\maketitle
\date
\address{\noindent
\kern -18pt
Rapha\"el Cerf,
DMA, \'Ecole Normale Supérieure, CNRS, PSL University, 75005 Paris.\\
LMO, Université Paris-Sud, CNRS, Université Paris–Saclay, 91405 Orsay.\\
\email{Raphael.Cerf@ens.fr}
\and
Carlo Mariconda (corresponding author), ORCID $0000-0002-8215-9394$,  Universit\`a
degli Studi di Padova,  Dipartimento di Matematica  ``Tullio Levi-Civita'',   Via Trieste 63, 35121 Padova, Italy; \email{carlo.mariconda@unipd.it}}
\begin{abstract}
\noindent
Let $L:\R\times \R\to [0, +\infty[\,\cup\{+\infty\}$ be a Borel function. 
We consider the problem 
\begin{equation}\tag{P}\min F(y)=\int_0^1L(y(t), y'(t))\,dt: y(0)=0,\, y\in W^{1,1}([0,1]).\end{equation}
We give an example of a real valued Lagrangian $L$ for which the Lavrentiev phenomenon occurs. We state a condition, involving only the behavior of $L$ on the graph of two functions, that ensures the non-occurrence of the phenomenon.  Our criterium weakens substantially the well-known condition, that $L$ is bounded on bounded sets. 
\end{abstract}
\section{Introduction}
Consider the basic problem of the Calculus of Variations that consists on minimizing the autonomous integral functional
\[F(y)=\int_0^1L(t,y(t), y'(t))\,dt\]
among the absolutely continuous functions on $[0,1]$ that possibly satisfy some end-point conditions.
Here $L:[0,1]\times \R^n\times\R^n\to [0, +\infty[\cup\{+\infty\}$ is a Borel function.
We are concerned with the question of avoiding the Lavrentiev phenomenon, namely the unpleasant fact that the infimum of $F$ among the absolutely continuous functions is strictly less than the one among the Lipschitz functions that share the same end-point conditions.
The occurrence of this phenomenon implies the failure of  classical numerical analysis methods, e.g., finite elements, if one wishes to compute the infimum of $F$, and represents a discontinuity of $F$ with respect to strong convergence in $W^{1,1}(I,\R^n)$.

Lavrentiev's phenomenon is considered  among experts a matter of non-au\-to\-no\-mous Lagrangians, i.e., depending explicitly on the time variable.   
On one side, a famous example by Manià exhibits the phenomenon when $$L(t,y, v)\,=\,(y^3-t)^2v^6$$ among the functions $y:[0,1]\to \R$ that satisfy the end-point conditions $y(0)=0, y(1)=1$ (see \cite[\S 4.3]{GBH}).  A more refined construction by Ball and Mizel \cite{BallM} shows that it may even occur when the Lagrangian is  a polynomial in $(t,y,v)$ that satisfies Tonelli's existence conditions (namely superlinearity and convexity in the last variable).
On the other side, a celebrated result by Alberti and Serra Cassano \cite[Theorem 2.4]{ASC} asserts that non pathological autonomous Lagrangians do not exhibit the phenomenon. More precisely, 
there is no Lavrentiev phenomenon if
\begin{equation}\tag{B}\forall K>0\quad \exists r_K>0\quad L(y,v)\text{ is bounded on }[-K,K]^n\times [-r_K,r_K]^n\,.\end{equation}
Notice that Condition (B) forces $L$ to be finite on the union $\bigcup_{K>0}
[-K,K]^n\times [-r_K,r_K]^n$ and in particular on $\R^n\times {\{0\}}$.

Actually, it has now become clear that the phenomenon is { also} strictly related to the presence and the number of end-point constraints. For instance, as shown in \cite{GBH}, Manià's example does not exhibit any more the phenomenon if one considers just the end-point condition $y(1)=1$.  Moreover, it was  pointed out in \cite{CM5}  that Condition (B) in \cite[Theorem 2.4]{ASC} is a sufficient condition for the non-occurrence of the phenomenon when one considers just one end-point condition, but not {anymore} in the
case of two end-point conditions.
In fact, Alberti provided an example (see \cite[Example 3.5]{CM5}) showing an autonomous  Lagrangian with values either 0 or $+\infty$,   { satisfying (B) such that} the  functional $F$  takes the value $+\infty$ on every Lipschitz function satisfying $y(0)=0, y(1)=1$.

Regarding the Lavrentiev phenomenon, the difference between one and two end-point conditions was recently better understood {(see \cite{CM5})}. 
{It seems} to us of interest to study more thoroughly the conditions that provide the non-occurrence of the phenomenon for problems with one end-point condition. As mentioned above, Condition (B) of \cite[Theorem 2.4]{ASC} does not take into account the { geometry of the} effective domain $\operatorname{Dom}(L)$ of the Lagrangian (i.e., the set where it is finite). We wonder how sharp condition (B) is.
Can it be weakened to an assumption involving { just the subsets of the effective domain of $L$? The effort of finding such a condition was carried out in \cite{CM7} under an additional  convexity hypothesis {along the radii from the origin} on the last variable of the Lagrangian: it is enough in that case that for each $K>0$, there is $r_K>0$ such that $L(y,v)$ is bounded on $\left([-K,K]^n\times [-r_K,r_K]^n\right)\cap \operatorname{Dom}(L)$.} 

In this paper, we consider the case where the Lagrangian $L=L(y,v)$ is autonomous, $n=1$, with the  initial condition $y(0)=0$ and free end-point condition.
We first exhibit a finite, autonomous Lagrangian that violates (B), for which the Lavrentiev phenomenon occurs with just one end-point condition. We introduce the 
following Condition (R), weaker than (B), that ensures, {with no need of any other additional hypothesis,} the non-occurrence of the phenomenon:
\begin{CondR}
 There exist two locally Lipschitz functions $\rho^-, \rho^+$ defined on $\mathbb R$
 such that:
 $$\forall z\in \R\qquad \rho^-(z)<0\,,\quad \rho^+(z)>0\,,$$
 and for every bounded interval $J$ of $\R$,
 $$\sup_{z\in J}\,L\big(z,\rho^-(z)\big)\,<\,+\infty\,,\quad
  \sup_{z\in J}\,L\big(z,\rho^+(z)\big)\,<\,+\infty\,.$$
\end{CondR}
\noindent
The Condition (R) is fulfilled when (B) of \cite[Theorem 2.4]{ASC} holds. {Condition (R) has the advantage to require the boundedness of the Lagrangian just on some one-dimensional subsets of its effective domain, without imposing, as (B) does, that $\operatorname{Dom}(L)$ contains the union of two-dimensional rectangles.}
\section{The functional, the gap and the phenomenon}
\subsection{The functional}
{For $1\le p\le +\infty$ we will denote by $W^{1,p}(I)$  the Sobolev space of absolutely continuous functions $y:I\to\R$ such that $y'\in L^p(I)$; $\Lip(I)=W^{1, \infty}(I)$ is the space of  Lipschitz functions on $I$.}

We consider an autonomous Borel Lagrangian $L:\R\times\R\to [0, {+\infty]}$ with non--negative values,
possibly infinite.
We denote by $I$ the unit interval $[0,1]$ and we define
\[\forall y\in W^{1,1}(I)\qquad F(y)\,=\,\int_0^1L(y(t), y'(t))\,dt.\]
We consider the end--point condition $y(0)=0$ and
the problem {\rm (P)}:
\begin{equation}\tag{P}
\min F(y)\,=\,\int_0^1L(y(t), y'(t))\,dt\,,\quad 
y\in W^{1,1}(I)\,,\quad 
y(0)=0\,.\end{equation}
\subsection{The Lavrentiev gap and phenomenon}
{We focus our attention on the Lavrentiev gap and phenomenon problems with the prescribed boundary condition  $y(0)=0$. The following definition can be easily adapted to other kinds of boundary conditions.}
\begin{definition}[No gap]
Let $p\geq 1$.
Let $y\in W^{1,p}(I)$ be such that $F(y)<+\infty$. 
We say that the {\em Lavrentiev gap} does not occur at $y$ for {\rm (P)} if
there exists a sequence $\left(y_{{h}}\right)_{{h}\in\N}$ of functions 
satisfying:
 \begin{enumerate}
\item for each $h\in \N$, the function $y_h$ is Lipschitz and $y_h(0)=0$;
 \item $\displaystyle\lim_{{
 h}\to +\infty}F(y_{h})= F(y)$ (approximation in {\em energy});
 \item $y_{h}\to y$ in $W^{1,p}(I)$ (approximation in {\em norm}).
\end{enumerate}
\end{definition}
\noindent
We denote by $\Lip(I)$ the space of the Lipschitz functions defined on $I$ with values in $\R$.
\begin{definition}[No phenomenon]
We say that there is no {\em Lavrentiev phenomenon} for ${\rm (P)}$  if
\[\inf ({\rm P})\,=\,\inf\,\big\{\,F(y):\, y\in\Lip(I)\,,y(0)=0\,\big\}\,.\]
\end{definition}
{
\begin{remark}
Clearly, the Lavrentiev phenomenon does not occur for ${\rm (P)}$ once there is no Lavrentiev gap for every $y\in W^{1,{p}}(I)$ such that $F(y)<+\infty$. 
Some others kinds of gap are of interest, namely between $W^{1,p}(I)$ and $W^{1,q}(I)$ for some $p<q<+\infty$. They are not considered here since our main result  directly prevents the occurrence of the gap between $W^{1,1}(I)$ and $W^{1,+\infty}(I)$: the technique of the proof of Theorem \ref{thm:Lav11} requires the functions $\rho^+$ and $\rho^-$ to be locally Lipschitz.
\end{remark}}
\section{Occurrence of the Lavrentiev gap for {\rm (P)}}
{The Lavrentiev phenomenon is often considered a pathology related to non-auto\-no\-mous Lagrangians.}
However, the phenomenon may also occur in the { autonomous} case: an example due to Alberti (see \cite[Example 3.5]{CM5}) exhibits an autonomous  Lagrangian $L$ that takes the value $+\infty$, { satisfies (B), yet}  the Lavrentiev phenomenon occurs for
a problem with two end-point conditions.
{When Condition (B) fails, the phenomenon may occur in the autonomous case, when one considers just one end-point condition. Consider}  
\[L(y, v)\,=\,\begin{cases}\left(v^2-\dfrac1{4y^2}\right)^2v^2&\text{ if } y\not=0\,,\\
\qquad 1&\text{ if }y=0\,.
\end{cases}\]
\begin{figure}[h!]
\begin{center}
\includegraphics[width=0.6\textwidth]{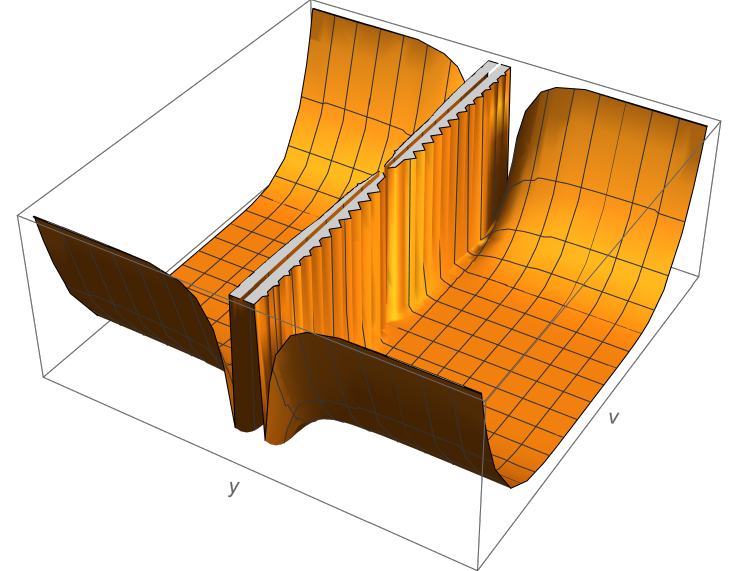}
\caption{{\small The graph of $L$ in Proposition \ref {pgap}.}}
\end{center}
\end{figure}
The Lagrangian $L$ is a Borel non--negative map. 
{ Let $y_*(s)=\sqrt{s}, s\in [0,1]$.
Notice that 
$(L, y_*)$ violates Condition (B) in \cite{ASC}.
Indeed if $v\not=0$, then
\[\lim_{y\to 0}L(y, v)\,=\,+\infty\,,\]
and this implies that, for every  $r>0$, $L$ is unbounded on $y_*(I)\times [-r,r]$.}
\begin{proposition}\label{pgap}
The function $y_*$ is a minimizer for the problem {\rm (P)} associated to $L$ and moreover $F(y_*)=0$. The gap occurs at $y_*$: more precisely,
for any $y\in\Lip([0,1])$ such that $y(0)=0$ { we have $F(y)\ge 1$}.
\end{proposition}
\begin{proof}
We have
$$\forall t>0\qquad y_*'(t)\,=\,\frac{1}{2\sqrt t}\,=\,\dfrac1{2y_*(t)}$$ so that
$L(y_*(t), y_*'(t))=0$ for $t$ in $]0,1]$. Moreover, 
we have $F(y)\ge 0$
for all admissible trajectory $y\in W^{1,1}([0,1])$, therefore
\[ F(y_*)=0=\min {\rm (P)}\,.\]
\noindent
Next, let $y\in\Lip([0,1])$ {with $y(0)=0$ and assume that {$F(y)<1$}.} {Notice that, since $F(0)=1$, then  $y$ is not identically equal to 0.}  The set $\{t\in [0,1]:\, y(t)\not=0\}$ being open and 
non-empty, it is a countable or finite  union  of non-empty open subintervals of $[0,1]$. Hence, there are $0\le a<b\le 1$ such that
\[y(a)=0,\quad y(b)\not=0,\quad \forall t\in]a,b[\quad y(t)\not=0\,.\]
Let $c\in ]a,b[$.
We have
\begin{equation}\label{hhuhuhi}\begin{aligned}
F(y)&\ge \int_c^bL(y(t), y'(t))\,dt\\&=\int_{c}^b\left(y'(t)^6-\dfrac12\dfrac{y'(t)^4}{y(t)^2}
+\frac{1}{16}\dfrac{y'(t)^2}{y(t)^4}\right)\,dt\\
&\ge -\dfrac12\int_c^b\dfrac{y'(t)^4}{y(t)^2}\,dt
+\frac{1}{16}\int_{c}^b\dfrac{y'(t)^2}{y(t)^4}\,dt\\
&\ge
-\dfrac12\|y'\|_{\infty}^2\int_{c}^b\dfrac{y'(t)^2}{y(t)^2}\,dt
+\frac{1}{16}\int_c^b\dfrac{y'(t)^2}{y(t)^4}\,dt\,.
\end{aligned}\end{equation}
Since
$y(a)=0$ and $y$ is continuous at $a$, then $y(t)\to 0$ as $t\to a$, so that
there exists {$d\in ]a,b[$}
such that
\begin{equation}\label{tag:ineqyy'}
\forall t\in ]a,d]\qquad
\dfrac12\|y'\|_{\infty}^2\dfrac{y'(t)^2}{y^2(t)}\,
\leq\,
\frac{1}{32}\dfrac{y'(t)^2}{y(t)^4}\,.
\end{equation}
{Now, in \eqref{hhuhuhi}, fix $c\in ]a,d[$}.
Integrating {both terms of \eqref{tag:ineqyy'}} over $[c,d]$, we obtain
\begin{equation}
\label{ert}
-\dfrac12\|y'\|_{\infty}^2\int_{c}^d\dfrac{y'(t)^2}{y^2(t)}\,dt\,
+\,
\frac{1}{16}\int_c^d\dfrac{y'(t)^2}{y(t)^4}\,dt\,\geq\,
\frac{1}{32}\int_c^d\dfrac{y'(t)^2}{y(t)^4}\,dt
\,.
\end{equation}
Inequalities~\eqref{hhuhuhi} and~\eqref{ert} together yield
\begin{equation}\label{jjuhuhi}
F(y)\,\ge\,
-\dfrac12\|y'\|_{\infty}^2\int_{d}^b\dfrac{y'(t)^2}{y^2(t)}\,dt
+\frac{1}{16}\int_d^b\dfrac{y'(t)^2}{y(t)^4}\,dt\,+\,
\frac{1}{32}\int_c^d\dfrac{y'(t)^2}{y(t)^4}\,dt
\,.
\end{equation}
Notice that we use the lower bound~\eqref{ert} only for the integral over $[c,d]$.
Jensen's inequality yields
\begin{equation}\label{tuhi}\begin{aligned}
\int_{c}^d\dfrac{y'(t)^2}{y(t)^4}\,dt&\,\ge\, \dfrac1{d-c}\left(\int_c^d
\dfrac{y'(t)}{y(t)^2}\,dt\right)^2\\
&\,=\,\dfrac1{d-c}\left(
\dfrac1{y(c)}-\dfrac1{y(d)}
\right)^2\,.
\end{aligned}\end{equation}
Since
$y(a)=0$ and $y$ is continuous at $a$, we deduce from~\eqref{tuhi} that
\begin{equation}\label{tag:bbhi}
\lim_{c\to a}\,\int_{c}^d\dfrac{y'(t)^2}{y(t)^4}\,dt\,=\,+\infty\,.
\end{equation}
Keeping $d$ fixed and
taking the limit in
\eqref{jjuhuhi}
as $c$ goes to $a$, we conclude that
$F(y)=+\infty$, { contradicting the initial assumption that $F(y)<1$}.
\end{proof}

\noindent
Proposition~\ref{pgap} readily implies that
the Lavrentiev phenomenon occurs for the problem
with one end-point condition given by
\begin{equation}
\min \int_0^1L(y(t), y'(t))\,dt\,,\quad 
y\in W^{1,1}(I)\,,\quad y(0)=0\,.
\end{equation}
The recent works \cite{CM7, CM5} 
have shown that 
the Lavrentiev phenomenon might be very sensitive to the number of { prescribed} end-point conditions.
In fact, the same argument as above shows that
the Lavrentiev phenomenon occurs for the problem
with two end-point conditions given by
\begin{equation}
	\min \int_{0}^1L(y(t), y'(t))\,dt\,,\quad 
	y\in W^{1,1}([0,1])\,,\quad y(0)=-1\,,y(1)=1\,.
\end{equation}
However, in the above example, if we keep only one of the two end-point conditions, 
the Lavrentiev phenomenon disappears! To see this, we proceed as in \cite[\S 4.3]{GBH}.
\section{Non-Occurrence of the Lavrentiev Gap and phenomenon for {\rm (P)}}\label{sect:P1X}
\subsection{Non-occurrence of the gap}
Let $L:\R\times\R\to [0, {+\infty]}$ 
be an autonomous Borel Lagrangian 
with non--negative values, possibly infinite.
Let $p\geq 1$.
For a given $y\in W^{1,p}(I)$, we consider the  following condition.
\begin{CondRy} There exist two Lipschitz functions $\rho^-, \rho^+$ defined on $y(I)$
 such that 
 $$\forall z\in y(I)\qquad \rho^-(z)<0\,,\quad \rho^+(z)>0\,,$$
 $$\sup_{z\in y(I)}\, L\big(z,\rho^-(z)\big)\,<\,+\infty\,,\quad
  \sup_{z\in y(I)}\,L\big(z,\rho^+(z)\big)\,<\,+\infty\,.$$
\end{CondRy}
\begin{theorem}[\textbf{Non-occurrence of the Lavrentiev  gap}]\label{thm:Lav11} 
Let $y\in W^{1,p}([0, 1])$ be such that $F(y)<+\infty$. Assume that 
 $y$ satisfies Condition (R$_y$).
Then there is {no Lavrentiev gap} for ${\rm (P)}$ at $y$.
\end{theorem}
{\begin{remark}
Condition (R$_y$) weakens Assumption (B) formulated in \cite{ASC}. Indeed if $L$ is bounded on $y(I)\times [-r, r]$ for some $r>0$, then (R$_y$) is satisfied with $\rho^-=-r$ and $\rho^+=r$. 
\end{remark}}
The strategy of the proof is the same as the proof of Alberti and Serra Cassano \cite{ASC}.
Yet it differs at some specific points and it requires also a different construction for
the approximating function. In order to facilitate the reading, we have chosen to write
the full proof. 
Another reason is that, as we work with real valued functions, some arguments become simpler than in the $n$--dimensional case.
For convenience, we restate {a} general lemma of integration
theory that were proved in \cite{ASC}.
\begin{lemma}\label{seca1}\cite[Lemma 2.6]{ASC}
Let $g:I\to[0,+\infty]$ be a Lebesgue measurable
function and let $B_h$ be a sequence of measurable subsets
of $I$ such that $|I\setminus B_h|\to 0$ as $h\to\infty$.
Then
$$\lim_{h\to+\infty}
\int_{B_h}g\,dt\,=\,
\int_{I}g\,dt\,.$$
\end{lemma}
{We will need the following variation of \cite[Lemma 2.7]{ASC}. { Compared with the original version we do not need here that the functions $\varphi_h$ are Lipschitz nor that $\varphi_h'$ are bounded from below on $I$, and the conclusion is weaker.}
\begin{lemma}\label{seca2}
Let $\varphi_h:I\to\R$, $h\in\mathbb N$, be a sequence of absolutely continuous functions with $\varphi_h'>0$ a.e. on $I$ and $(E_h)_{h\in\N}$ be a sequence of measurable sets of $I$ such that:
\begin{enumerate}
\item 
 For some $c>0$, $\varphi'_h\geq c$ a.e. in $E_h$;
 \item 
$\varphi_h(t)\to t$
as $h\to+\infty$ for every $t\in I$.
\end{enumerate}
Then, for every $f\in L^p(\R)$, 
\[\lim_{h\to +\infty}\int_{\varphi_h^{-1}(I)\cap E_h}|f-f(\varphi_h)|^p\,d\tau=0.\]
\end{lemma}
\begin{proof} Let $(f_k)_{k\in\N}$ be a sequence of smooth functions such that $f_k\to f$ in $L^p(I)$. We write
\begin{multline}\label{tag:badlemma}\|f-f(\varphi_h)\|_{L^p(\varphi_h^{-1}(I)\cap E_h)}\le\\
\|f-f_k\|_{L^p(I)}+\|f_k-f_k(\varphi_h)\|_{L^p(I)}+\|f_k(\varphi_h)-f(\varphi_h)\|_{L^p(\varphi_h^{-1}(I)\cap E_h)}.\end{multline}
Of course, $\|f-f_k\|_{L^p(I)}\to 0 $ in $L^p(I)$ and 
the
dominated convergence theorem, for each $k$ fixed, implies that
$\|f_k-f_k(\varphi_h)\|_{L^p(I)}\to 0$ as $h\to +\infty$. 
It remains to study the convergence of $\|f_k(\varphi_h)-f(\varphi_h)\|_{L^p(\varphi_h^{-1}(I)\cap E_h)}$.
For each $h\in\mathbb N$ let $\psi_h:\varphi_h(I)\to I$ be the inverse of $\varphi_h$; notice that by {Banach-Zarecki}'s theorem $\psi_h$ is absolutely continuous (see \cite{Bruckner, Natanson}).
The change of variables $t=\varphi_h(\tau)$ yields (see, for instance, \cite{Serrin}),
{
\[\begin{aligned}
\int_{\varphi_h^{-1}(I)\cap E_h}|f_k(\varphi_h)-f(\varphi_h)|^p\,d\tau&\leq
\frac{1}{c}\int_{\varphi_h^{-1}(I)\cap E_h}|f_k(\varphi_h)-f(\varphi_h)|^p\varphi'_h\,d\tau\\
&\leq
\frac{1}{c}\int_{\varphi_h^{-1}(I)}|f_k(\varphi_h)-f(\varphi_h)|^p\varphi'_h\,d\tau\\
&= 
\dfrac1c\int_{I}|f_k(t)-f(t)|^p\,dt\\&=
\dfrac1c
\|f_k(t)-f(t)\|^p_{L^p(I)}\,dt\to 0
\end{aligned}\]
}
as $k\to +\infty$.
The conclusion follows.   
\end{proof}
}
\noindent
We will use several times the fact that, since $\rho^+, \rho^-$ are continuous functions on $I$,
there are positive constants $\rho_{min}>0,
\rho_{max}>0$ 
such that 
\[
\forall x\in I\qquad
\min\{\rho^+(x), -\rho^-(x)\}\ge \rho_{min}, \quad  \max\{\rho^+(x), -\rho^-(x)\}\le \rho_{max}
\,.\]
\begin{proof}
	[{Proof of Theorem \ref{thm:Lav11}}]
We assume that $y$ is not constant, otherwise the conclusion is trivial.
We start by applying a classical result which is a consequence
of Lusin's theorem 
(see for instance Theorem~$3.10$ in \cite{Z}).
 For every $h\in\mathbb N$,
	there are a Lipschitz function $u_h:I\to \R$ and an open subset $A_h$ of $I$ such that, { denoting by $|A|$ the measure of a subset $A$ of $I$}:
       \begin{enumerate}[leftmargin=*]
   \item  $|A_h|\leq 1/h$,
   \item  $u_h=y,\, u_h'=y'$ in $I\setminus A_h$,
   \item $u_h(0)=y(0)$, $u_h(1)=y(1)$,
   \item  $u_h$ is affine in each connected component of $A_h$,
   \item $A_h$ is a countable
union of disjoint open intervals $I_{h,k}$, 
$k\in {J_h}\subset\N$. 
 \end{enumerate}
The set $A_h$ is somehow the bad set where the function $y$
 might behave badly in the sense that its derivative might be unbounded { on $A_h$.}\\
{We claim that it is not restrictive to 
assume that $u_h'$ does not vanish on $A_h$.
Indeed, assume the contrary. Each $A_h$ is a union of disjoint open intervals $\left(I_{h,k}\right)_{ k\in J_h}$, and we may  modify the sets $A_h$  as follows:}
 \begin{itemize}[leftmargin=*]
     \item We first remove from $A_h$ the intervals where $y$ is itself constant on which, as a byproduct, $y=u_h$;
     \item On every other subinterval $I_{h,k}=(a_{h,k},b_{h,k})$, ${k\in J'_h\subset J_h}$ where $u_h$ is constant but $y$ is not, we choose $c_{h,k}\in I_{h,k}$ such that $y(c_{h,k})\not=y(a_{h,k})$.  On {\[I_{h,k}\setminus\{c_{h,k}\}=(a_{h,k}, c_{h,k})\cup (c_{h,k}, b_{h,k})\]} we define $\tilde {u}_h$ to be affine in $(a_{h,k}, c_{h,k})$ joining $y(a_{h,k})$ to $y(c_{h,k})$ and   affine in $(c_{h,k}, b_{h,k})$ joining $y(c_{h,k})$ to $y(b_{h,k})$. {Namely
     \[\forall\tau\in [a_{h,k}, c_{h,k}]\, \quad {u}_h(\tau)=y(a_{h,k})+\dfrac{\tau-a_{h,k}}{c_{h,k}-a_{h,k}}(y(c_{h,k})-y(a_{h,k})),\]
     \[\forall\tau\in [c_{h,k}, b_{h,k}]\, \quad {u}_h(\tau)=y(c_{h,k})+\dfrac{\tau-c_{h,k}}{b_{h,k}-c_{h,k}}(y(b_{h,k})-y(c_{h,k})).\]} 
     \end{itemize}
{ We then set $$\tilde {A}_h=A_h\setminus \bigcup_{ k\in J'_h}\{c_{h,k}\}\,.$$
     Clearly $|A_h|=|\tilde {A}_h|$ and $(\tilde{u}_h, \tilde {A}_h)$ satisfy properties 1-5,  proving the claim: we thus assume henceforth that  $u_h'$ does not vanish on $A_h$.}\\
Notice that, since $u_h$ is affine on every interval $I_{h,k}$ and $u_h,y$ are equal at the extremities of $I_{h,k}$, then
 \[\begin{aligned}\label{boundb}
 \int_{A_h}|u_h'|\,d\tau&\,=\sum_{k\in J_h}\int_{I_{h,k}}|u_h'|\,d\tau
 \,=\,\sum_{k\in J_h}\left|\int_{I_{h,k}}u_h'\,d\tau\right|\\
 &\,=\sum_{k\in J_h}\left|\int_{I_{h,k}}y'\,d\tau\right|\\
 &\,\le \sum_{k\in J_h}\int_{I_{h,k}}|y'|\,d\tau=\,\int_{A_h}|y'|\,d\tau\to 0\quad\text{as}\quad h\to +\infty.
 \end{aligned}\]
The first problem with the function $u_h$ is that we might have $F(u_h)=+\infty$,
indeed, the integral of $L{(u_h, u_h')}$ over the intervals $I_{h,k}$ might very well
be infinite. We shall take advantage of the functions $\rho^-,\rho^+$ to
replace the portions of the function $u_h$ over $A_h$ by a function $v_h$
having a finite energy on $A_h$, which in addition
tends to 0 as $h\to +\infty$. 
We define a function $\rho_h$ {  on
$A_h$  by setting}
\begin{equation}
	\label{defrho}
\forall \tau\in A_h\qquad
\rho_h(\tau)\,=\,\begin{cases}
	\phantom{+}\rho^+(u_h(\tau))&\text{ if } u_h'(\tau)> 0\,,\\
-\rho^-(u_h(\tau))&\text{ if }u_h'(\tau)<0\,.\end{cases}\end{equation}
{ Notice that $\rho_h$ is positive and continuous.}
In order to perform an adequate change of variable,
we define next a function $\varphi_h\in\ W^{1,1}(I)$ by
setting $\varphi_h(0)=0$ and
\begin{equation}
	\label{defphip}
	 \varphi_h'({\tau})=
 \begin{cases}\quad 1&\text{ if }\tau\in I\setminus A_h\,,\\
\dfrac{|u_h'(\tau)|}{\rho_h(\tau)} &\text{ if } \tau\in A_h\,.\end{cases}
\end{equation}
Using inequality~\eqref{boundb}, we have
  \[\begin{aligned}\label{tag:zfestimate}|
  \varphi_h(A_h)|&\,=\,\int_{\varphi_h(A_h)}1\,d{t}\,=\,
  \int_{A_h}\varphi_h'(\tau)\,d\tau\,=\,
  \int_{A_h}\frac{|u'_h(\tau)|}{\rho_h(\tau)}\,d\tau\\
     &\,\le\, \dfrac1{\rho_{min}}\int_{A_h}|u_h'(\tau)|\,d\tau
     \,\to\, 0 \quad\text{as}\quad h\to +\infty\,.
     \end{aligned}\]
     Next, we claim that
the function $\varphi_h$  converges uniformly towards the identity
 map on $I$. Indeed, $\varphi_h(0)=0$ and moreover, {from \eqref{tag:zfestimate} and the fact that $|A_h|\to 0$,} { for all $t\in I$ we have}
\begin{equation}\label{idenca}\begin{aligned}
{|\varphi_h(t)-t|}&\le \int_I|\varphi_h'-1|\,d\tau\,\\\\&\leq
\int_{A_h}\left(\varphi_h'(\tau)+1\right)\,d\tau\\
  \,&=\, |\varphi_h(A_h)|+|A_h|\,\to\, 0\phantom{\int}
  \quad\text{as}\quad h\to +\infty\,.\end{aligned}\end{equation}
In particular, we have 
 $$|\varphi_h(I)|\,=\,\displaystyle\int_I\varphi_h'(\tau)\,d\tau\,=\,
 \varphi_h(1){-\varphi_h(0)}\,\to \,|I|=1\quad \text{as}\quad h\to +\infty\,.$$
Setting $T_h=\varphi_h(1)$, we thus have $T_h\to 1$ as $h\to +\infty$.
However we don't know whether $T_h$ is smaller or larger than~$1$ and this will create
some trouble later on.
The derivative $\varphi'_h$ is strictly positive  on $I$, therefore
$\varphi_h$ is  strictly increasing and it is a one to one map from $[0,1]$ onto $[0,T_h]$.  
Its inverse  $\psi_h:[0, T_h]\to [0,1]$ is continuous, strictly increasing. { Since $\varphi'_h>0$ a.e. on $I$ then $\psi_h$ is absolutely continuous (see \cite[Ch. IX]{Natanson})} and its  derivative is given by
\begin{equation}
	\label{ihi}{\forall t\in [0, T_h]\qquad}
	\psi_h'(t)\,=\,\frac{1}{\varphi'_h\big(\psi_h(t)\big)}\,.
\end{equation}
Using the expression of $\varphi'_h$ given in \eqref{defphip}, we obtain
\begin{equation}
	\label{invphi}
\psi_h'(t)\,=\,\begin{cases}
\qquad 1&\text{ if }t\in [0, T_h]\setminus \varphi_h(A_h)\,,\\
\dfrac{\rho_h(\psi_h(t))}{|u_h'(\psi_h(t))|} &\text{ if } t\in \varphi_h(A_h)\,.\end{cases}
\end{equation}
The change of variable $\tau=\psi_h(t)$ gives
\[\begin{aligned}{1=}|I|&\,\ge\, |\psi_h(I\cap [0, T_h])|
	\,=\,\int_{\psi_h(I\cap [0, T_h])}{1}\,d\tau
	\,=\,\int_{I\cap [0, T_h]}\psi'_h(t)\,dt\\
&\,\ge\, \int_{{(I\cap [0, T_h])}\setminus \varphi_h(A_h)} 1\,dt
\,\ge\, \min\{1, T_h\}-|\varphi_h(A_h)|\,,\end{aligned}\]
so that
\begin{equation}
	\label{psih}
	\lim_{h\to+\infty}\,
	|\psi_h(I\cap [0, T_h])|\,=\, 1\,.
\end{equation}
We define a new function $v_h$ by setting
\[\forall t\in [0, T_h]\qquad v_h(t)=u_h(\psi_h(t)).\]
The function $v_h$, being the composition of the Lipschitz function $u_h$ 
with the absolutely continuous function $\psi_h$, is absolutely continuous with derivative given by
\[\forall t\in [0, T_h]\qquad 
v_h'(t)\,=\,u_h'(\psi_h(t))\,\psi_h'(t)\,.\]
For $t\in [0, T_h]\setminus \varphi_h(A_h)$, we have
$\psi'_h(t)=1$ and
$\psi_h(t)\not\in A_h$, whence
\begin{equation}
	\label{avr}
	u_h'(\psi_h(t))\,\psi'_h(t)
\,=\,u_h'(\psi_h(t))
\,=\,y'(\psi_h(t))\end{equation} and thus
\begin{equation}
	\label{avq}
	\forall t\in 
[0, T_h]\setminus \varphi_h(A_h)
	\qquad 
v_h'(t)\,=\,y'(\psi_h(t))\,.\end{equation}
Let next
$t\in \varphi_h(A_h)$. In this case, we have
\begin{equation}\label{nhu}
\begin{aligned}
v_h'(t)&\,=\,
	     \dfrac{u'_h(\psi_h(t))}{\varphi_h'(\psi_h(t))}
	    \,=\,
	     \dfrac{u'_h(\psi_h(t))}{|u_h'(\psi_h(t))|}\rho_h(\psi_h(t))\\
	    &\,=\,
	     \sgn(u'_h(\psi_h(t)))\,\rho_h(\psi_h(t))\,,\phantom{\dfrac{u'_h(\psi_h(t))}{|u_h'(\psi_h(t))|}\rho_h(\psi_h(t))}
	     \end{aligned}
\end{equation}
where $\sgn$ is the classical sign function given by
$$\sgn(x)\,=\,\begin{cases}+1&\text{ if }x> 0\,,\\-1&\text{ if } x<0\,.\end{cases}$$
Recalling the definition~\eqref{defrho} of $\rho_h$,
we see that
\begin{equation}
\forall t\in \varphi_h(A_h) 
\qquad
v'_h(t)\,=\,
\begin{cases}
	\rho^+(u_h(\psi_h(t)))&\text{ if}\quad u_h'(\psi_h(t))> 0\,,\\
\rho^-(u_h(\psi_h(t)))&\text{ if}\quad u_h'(\psi_h(t))<0\,.\end{cases}
\end{equation}
This implies in particular that  
\begin{equation}\label{stimvh}
\forall t\in \varphi_h(A_h)\qquad
|v_h'(t)|\,\le\, \rho_{max}\,.\end{equation}
From formula~\eqref{avr}, we see that $v'_h$ is also bounded on
$[0, T_h]\setminus \varphi_h(A_h)$, since
the function $u_h$ is Lipschitz.
We conclude that $v_h$ is Lipschitz on $[0, T_h]$.\\
{{\em Construction of the approximating sequence.}}
We finally { build  the Lipschitz function $w_h$ upon $v_h$}  { on $[0,1]$} which { is the suitable approximation of} $y$ in energy and in the space $W^{1,{p}}(I)$.
{Two cases may occur:\\
a) If $T_h\ge 1$, then we define $w_h$ to be the restriction of $v_h$ to $[0,1]$.\\
b) If $T_h< 1$, then we shall extend $v_h$ from $[0, T_h]$ to $[0,1]$ as we explain next.}
We define$$\alpha\,=\,\min\,y(I)\,,\qquad
\beta\,=\,\max\,y(I)\,.$$
{ We first extend $\rho^-, \rho^+$ to Lipschitz functions $\widetilde\rho^-, \widetilde\rho^+$ on $\R$ defined by
\[\widetilde\rho^-(z)=\begin{cases}\rho^-(\alpha)&z<\alpha\\\rho^-(z)&z\in [\alpha, \beta]
\\ \rho^-(\beta)&z>\beta.\end{cases},\quad \widetilde\rho^+(z)=\begin{cases}\rho^+(\alpha)&z<\alpha\\\rho^+(z)&z\in [\alpha, \beta]
\\ \rho^+(\beta)&z>\beta.\end{cases}\]
We consider the differential equation
$$z'(t)\,=\,\widetilde\rho^+(z(t))$$
	and we denote by $z^+_1(t)$ the solution
	with initial condition $z(\tau_0)=y_0$, 
	where $\tau_0=T_h$ and $ y_0=v_h(T_h)$.
We set
	$$\tau_1\,=\,\inf\,\big\{\,t\geq \tau_0: z^+_1(t)= \beta
	\,\big\}\,.$$
	If $\tau_1<1$, then we set 
	$w_h(t)=z_1^+(t)$ on $[\tau_0, \tau_1]$. 
	We consider then the differential equation
$$z'(t)\,=\,\widetilde\rho^-(z(t))$$
	and we denote by $z^-_2(t)$ the solution
	with initial condition $z(\tau_1)=z^+_1(\tau_1)=\beta$.
We set
	$$\tau_2\,=\,\inf\,\big\{\,t\geq \tau_1: z^-_2(t)= \alpha
	\,\big\}\,.$$
	Notice that, since $-\rho_{max}\le \rho^-$, 
	the travelling speed to go from $\beta$ to $\alpha$
	is at most $\rho_{max}$ and thus
	\[\tau_2-\tau_1\ge \dfrac{\beta-\alpha}{\rho_{max}}.\]
	If $\tau_2<1$, then we extend $w_h(t)$ on $[\tau_1,\tau_2]$ by setting
	$w_h(t)=z_2^-(t)$ on this interval. 
	We iterate this construction. 
	Since at each stage $i\ge 1$, we have
	\[\tau_{i+1}-\tau_i\ge \dfrac{\beta-\alpha}{\rho_{max}},\]
	then 
	the process ends {  at the first index $m$ such that
 $\tau_m<1\le \tau_{m+1}$}, after a number $m$ of steps that is bounded by a number depending only on $\beta- \alpha$ and $1-T_h$.
	In fact, we have
	\begin{equation}\label{boundm}
		m\,\leq\,
	 \dfrac{\rho_{max}}{\beta-\alpha}(1-T_h)+1\,.
\end{equation}
In the last step	 we extend 
	$w_h(t)$ on $[\tau_m,1]$ by restricting the solution
	of the differential equation
	to this interval.  In what follows we set 
	$\tau_{m+1}=1$ for convenience.
To sum up, we have 
 \[\forall t\in [0,1]\quad w_h(t)\in[\alpha, \beta],\quad w_h'(t)\,=\,
	 \begin{cases}
		 \quad v'_h(t)&\text{ if }t\in [0,T_h]\,,\\
		 \rho^+(w_h(t))&\text{ if }t\in [\tau_i, \tau_{i+1}]\,,\quad i\text{ even}\,,\\
		 \rho^-(w_h(t))&\text{ if }t\in [\tau_i, \tau_{i+1}]\,,\quad i\text{ odd}\,.
 \end{cases}\]}

Notice that 
\begin{equation}
	\label{whm}
	\forall t\in
[T_h, 1] \qquad
	|w_h'(t)|\,\le \,\rho_{max}\,,
\end{equation}
hence the function $w_h$ is still Lipschitz.

We show next that $w_h$ converges to $y$ in $W^{1,p}(I)$.  
 We decompose the integral as the sum of three terms
 \begin{equation}\label{tag:allequations1}
  {\|w'_h-y'\|^p_{L^p(I)}}\,=\,
  \int_{I}
|w'_h(t)-y'(t)|^p\,dt\,=\,
P_{1,h}+P_{2,h}+P_{3,h}\,,\end{equation}
where, recalling that $w_h=v_h$ on $[0,T_h]$,
  \begin{align}
	  P_{1,h}\,&=\int_{{(I\cap [0, T_h])}\setminus\varphi_h(A_h)}
|v'_h(t)-y'(t)|^p\,dt\,,\\
	  P_{2,h}\,&=\int_{I\cap [0, T_h]\cap \varphi_h(A_h)}
|v'_h(t)-y'(t)|^p\,dt\,,\\
  P_{3,h}\,&=\,\int_{{[{\min\{T_h,1\}},1]}}|w'_h(t)-y'(t)|^p\,dt\,.
  \end{align}
We prove next that the three terms $P_{1,h}, P_{2,h}, P_{3,h}$ tend to 0 as $h\to +\infty$.
We have
 \begin{equation}P_{1,h} 
	 \,=\,
 \int_{{(I\cap [0, T_h])}\setminus\varphi_h(A_h)}|u'_h(\psi_h(t))\psi'_h(t)-y'(t)|^p\,dt\,.
\end{equation}
It follows from~\eqref{invphi} that $\psi'_h=1$ on
 ${(I\cap [0, T_h])}\setminus\varphi_h(A_h)$, therefore we can rewrite 
 $P_{1,h}$ as
 \begin{equation}P_{1,h} 
	 \,=\,
 \int_{{(I\cap [0, T_h])}\setminus\varphi_h(A_h)}{|u'_h(\psi_h(t))-y'(t)|^p}\psi'_h(t)\,dt
\,.\end{equation}
The change of variable $\tau=\psi_h(t)$ then yields 
 \[P_{1,h} \,=\,
 \int_{\psi_h(I\cap [0, T_h])\setminus A_h}|u'_h(\tau)-y'(\varphi_h(\tau))|^p\,d\tau\,.\]
Using the fact that 
$u'_h=y'$ on $I\setminus A_h$, we obtain
{\[P_{1,h} 
 \,=\,
 \int_{\psi_h(I\cap [0, T_h])\setminus A_h}\big|y'(\tau)-y'(\varphi_h(\tau))\big|^p\,d\tau.\]
 Notice that $\psi_h(I\cap [0, T_h])\setminus A_h=\varphi_h^{-1}(I)\cap E_h$, with $E_h=I\setminus A_h$ and that $\varphi_h'=1$ on $E_h$.
 By applying {L}emma~\ref{seca2} we obtain that $P_{1,h}\to 0$ as $h\to +\infty$.}
Concerning 
$P_{2,h}$, we
notice that
\[P_{2,h}\le 2^p\left(\int_{I\cap[0, T_h]\cap \varphi_h(A_h)} |v'_h(t)|^p\,dt+ \int_{I\cap[0, T_h]\cap \varphi_h(A_h)} |y'(t)|^p\,dt  \right).\]
It follows from \eqref{stimvh} and~\eqref{tag:zfestimate}
that 
\[\int_{I\cap [0, T_h]\cap \varphi_h(A_h)} |v'_h(t)|^p\,dt
	\,\le\,( \rho_{max})^p\,|\varphi_h(A_h)|\to 0\quad\text{as}\quad h\to +\infty,
\]
and the integrability of $|y'|^p$ {together with \eqref{tag:zfestimate}}  immediately gives,
\[\int_{I\cap [0, T_h]\cap \varphi_h(A_h)} |y'(t)|^p\,dt\to 0
	\quad\text{as}\quad h\to +\infty,
\]
so that $P_{2,h}\to 0$ as $h\to +\infty$.
Finally, 
we have
\[P_{3,h}\le 2^p\left(\int_{{ [{\min\{T_h,1\}},1]}}|w'_h(t)|^p\,dt+\int_{{ [{\min\{T_h,1\}},1]}}|y'(t)|^p\,dt\right).\]
As above, since $T_h\to 1$ as $h\to +\infty$, {the integrability of $|y'|^p$ gives}
\[\int_{{ [{\min\{T_h,1\}},1]}}|y'(t)|^p\,dt\to 0\quad
\text{as}\quad h\to +\infty.\]
Moreover, $|w_h'|\le \rho_{max}$ on $I\cap [T_h,1]$, so that
\[\int_{{ [{\min\{T_h,1\}},1]}}|w'_h(t)|^p\,dt
	\,\le\,( \rho_{max})^p
	{(1-\min\{T_h,1\})}\to 0\quad\text{as}\quad h\to +\infty.
\]
\noindent
We  show now that $F(w_h)$ converges to $F(y)$ 
as $h\to +\infty$. By definition, we have
$$F(w_h)\,=\,
\int_{I}L(w_h, w_h')\,dt 
\,.$$
We decompose the integral as the sum of three terms
 $$F(w_h)\,=\,
Q_{1,h}+Q_{2,h}+Q_{3,h}\,,
$$
where, recalling that $w_h=v_h$ on $[0,T_h]$,
  \begin{align}
	  Q_{1,h}\,&=\,
    {\int_{{(I\cap [0, T_h])}\setminus\varphi_h(A_h)}L(v_h, v_h')\,dt }
	  \,,\\
	  Q_{2,h}\,&=\,
{ \int_{I\cap [0, T_h]\cap \varphi_h(A_h)}L(v_h, v_h')\,dt}
	  \,,\\
	  Q_{3,h}\,&=\,
{ \int_{{[{\min\{T_h,1\}},1]}}L(w_h, w_h')\,dt}
	  \,.\\
  \end{align}
We prove next that  $Q_{1,h}$ converges towards $F(y)$
while $Q_{2,h}, Q_{3,h}$ tend to 0 as $h\to +\infty$.
From the definition of $v_h$, we have
\[Q_{1,h}\,=\,\int_{(I\cap [0, T_h])\setminus\varphi_h(A_h)}
L\big(u_h(\psi_h(s)), u'_h(\psi_h(s))\psi'_h(s)\big)\,ds\,.\]
It follows from~\eqref{invphi} that $\psi'_h=1$ on
 ${(I\cap [0, T_h])}\setminus\varphi_h(A_h)$, therefore we can rewrite 
 $Q_{1,h}$ as{
\[Q_{1,h}\,=\,\int_{{(I\cap [0, T_h])}\setminus\varphi_h(A_h)}
L\big(u_h(\psi_h(s)), u'_h(\psi_h(s))\big)\,ds\,.\]}
The change of variable $\tau=\psi_h(s)$ yields then
\[Q_{1,h}\,=\,
\int_{\psi_h(I\cap [0, T_h])\setminus A_h}\kern-41pt
L\big(u_h(\tau), u_h'(\tau)\big)\,d\tau\,.\]
Using the fact that $u_h=y$ and
$u'_h=y'$ on $I\setminus A_h$, we obtain
\[Q_{1,h}\,=\,
\int_{\psi_h(I\cap [0, T_h])\setminus A_h}\kern-41pt
L\big(y(\tau), y'(\tau)\big)\,d\tau\,.\]
Lemma~\ref{seca1} and
the estimate~\eqref{psih} allow to conclude that
\[Q_{1,h}
\quad\to\quad \int_IL\big(y(\tau), y'(\tau)\big)\,d\tau\quad\text{as}\quad h\to +\infty\,.\]
Thus we are done with $Q_{1,h}$. We deal next with $Q_{2,h}$.
The expression of the derivative of $v'_h$ on $I\cap [0, T_h]\cap \varphi_h(A_h)$
was computed in \eqref{nhu}, so we have
\begin{equation}\label{eqs}Q_{2,h}\,=\,\int_{I\cap [0, T_h]\cap \varphi_h(A_h)}L\Big(u_h(\psi_h(t)), 	     \sgn(u'_h(\psi_h(t)))\,\rho_h(\psi_h(t))
\Big)\,dt\,.
\end{equation}
The change of variable $\tau=\psi_h(t)$ gives, with the
help of the expression of $\psi'_h$ computed in~\eqref{invphi},
\begin{equation}\label{tag:Q2}\begin{aligned}Q_{2,h}
&\,=\,\int_{\psi_h(I\cap [0, T_h])\cap A_h}L\Big(u_h(\tau), 
\sgn(u'_h(\tau))
\,\rho_h(\tau)
\Big)\dfrac{|u_h'(\tau)|}{\rho_h(\tau)}\,d\tau\\
	&\,\le \dfrac1{\rho_{min}}
	\int_{A_h}L\Big(u_h(\tau), 
\sgn(u'_h(\tau))
\,\rho_h(\tau)
	\Big)|u_h'(\tau)|\,d\tau\,.
\end{aligned}\end{equation}
{ We thus obtain
$$
Q_{2,h}\,\leq\,
	\dfrac1{\rho_{min}}
	\Big(
  \sup_{z\in y(I)}\,L\big(z,\rho^-(z)\big)
  +
  \sup_{z\in y(I)}\,L\big(z,\rho^+(z)\big)
  \Big)
	\int_{A_h}|u'_h|\,d\tau
  \,
		$$ so that { Condition $(R_y)$ and \eqref{boundb} yield} that $Q_{2,h}\to 0$ as $h\to +\infty$.}
		\\
It remains to prove that 
$Q_{3,h}\to 0$ as $h\to +\infty$. 
If $T_h\geq 1$, then $Q_{3,h}=0$.
Let us examine the case where
$T_h<1$.
From the construction of the extension $w_h$ of $v_h$ on $[T_h,1]$, we have
\begin{multline}Q_{3,h}\,=\,
	\sum_{\genfrac{}{}{0pt}{}{0\leq i\leq m}{ i\text{ even}}}\int_{\tau_i}^{\tau_{i+1}}L\Big(w_h(t), \rho^{+}(w_h(t))\Big)\,dt\cr
	\,+\,
	\sum_{\genfrac{}{}{0pt}{}{0\leq i\leq m}{ i\text{ odd}}}\int_{\tau_i}^{\tau_{i+1}}L\Big(w_h(t), \rho^{-}(w_h(t))\Big)\,dt
\end{multline}
{ so that 
$$
	Q_{3,h}\,\leq\,
	\Big(
  \sup_{z\in y(I)}\,L\big(z,\rho^-(z)\big)
  +
  \sup_{z\in y(I)}\,L\big(z,\rho^+(z)\big)
  \Big)(1-T_h).
		$$
Since $T_h\to 1$ as $h\to+\infty$, we conclude 
{ with the help of Condition $(R_y)$} that
	$Q_{3,h}\to 0$   as $h\to +\infty$.
		}
The proof that $ F(w_h)\to F(y)$ as $h\to +\infty$ is now complete.
\end{proof}
\noindent
Inspecting the proof of Theorem \ref{thm:Lav11}, we see 
that we could replace the Lipschitz continuity assumption on the functions 
$\rho^+, \rho^-$ by the assumption that 
{
$$\sup_{z\in y(I)}\,L\big(z,0\big)\,<\,+\infty.$$}
Indeed, in the case where $T_h<1$, we could { then} extend $v_h$ on $[T_h, 1]$ with the help of a constant function on $[T_h, 1]$. 
\subsection{Non-occurrence of the phenomenon}
{ We give finally a condition ensuring the non--occurrence of the Lavrentiev phenomenon.}
\begin{CondR}
 There exist two locally Lipschitz functions $\rho^-, \rho^+$ defined on $\mathbb R$
 such that:
 $$\forall z\in \R\qquad \rho^-(z)<0\,,\quad \rho^+(z)>0\,,$$
 and for every bounded interval $J$ of $\R$,
 $$\sup_{z\in J}\,L\big(z,\rho^-(z)\big)\,<\,+\infty\,,\quad
  \sup_{z\in J}\,L\big(z,\rho^+(z)\big)\,<\,+\infty\,.$$
\end{CondR}
\begin{corollary}[\textbf{Non-occurrence of the Lavrentiev  phenomenon}]\label{coro:Lav11} 
Let $L:\R\times\R\to [0, {+\infty]}$ 
be an autonomous Borel Lagrangian 
with non--negative values, possibly infinite.
Suppose that $L$ satisfies Condition (R). 
Then the Lavrentiev phenomenon does not occur for {\rm (P)}.
\end{corollary}
\begin{proof}
Let ${ (y_k)_{k\in\N}}$ be a minimizing sequence for {\rm (P)} satisfying, for all $k\in\mathbb N$:
\begin{itemize}
    \item  $y_k\in W^{1,1}(I)$, $y_k(0)=0$;
    \item $F(y_k)\le\inf{\rm (P)}+\dfrac1{k+1}$.
\end{itemize}
Fix $k\in\mathbb N$. The condition (R) implies that 
the condition (R$_{y_k}$) holds as well. 
By {T}heorem~\ref{thm:Lav11}, there exists
$z_k\in\Lip([0,1])$  such that
$z_k(0)=0$ and 
$$F(z_k)\,\le\, F(y_k)+\dfrac1{k+1}\,.$$
Therefore ${ (z_k)_{k\in\N}}$ is a minimizing sequence of Lipschitz functions for {\rm (P)}, thus proving the claim.
\end{proof}

\bibliographystyle{plain}
\bibliography{Lavrentiev1dim2}
\end{document}